\documentclass[11pt,reqno]{amsart}
\usepackage[T1]{fontenc}
\usepackage{lmodern}
\usepackage{amsmath,amssymb,mathtools}
\usepackage{esint}
\usepackage[paper=a4paper,textwidth=6.15in,textheight=8.65in]{geometry}
\usepackage{microtype}
\usepackage[colorlinks=true,linkcolor=blue,citecolor=blue,urlcolor=blue]{hyperref}
\hypersetup{pdftitle={An L1 bound for Ricci curvature on Kahler manifolds},
pdfauthor={Wenshuai Jiang and Gang Liu},
pdfsubject={Kahler Ricci L1 proof; complex structure J; revised notation}}
\allowdisplaybreaks[2]
\numberwithin{equation}{section}
\newtheorem{theorem}{Theorem}[section]
\newtheorem{proposition}[theorem]{Proposition}
\newtheorem{lemma}[theorem]{Lemma}
\newtheorem{corollary}[theorem]{Corollary}
\theoremstyle{definition}
\newtheorem{definition}[theorem]{Definition}
\theoremstyle{remark}
\newtheorem{remark}[theorem]{Remark}
\newtheorem*{remark*}{Remark}
\newcommand{\Ric}{\operatorname{Ric}}

\newcommand{\Vol}{\operatorname{Vol}}
\newcommand{\Hess}{\operatorname{Hess}}
\newcommand{\GH}{\operatorname{GH}}
\newcommand{\Lip}{\operatorname{Lip}}
\newcommand{\N}{\mathcal N}
\newcommand{\C}{\mathcal C}
\newcommand{\W}{\mathcal W}
\newcommand{\V}{\mathcal V}
\newcommand{\dd}{\,d}
\newcommand{\eps}{\varepsilon}
\newcommand{\fintt}{\fint}
\title[An \(L^1\) bound for Ricci curvature]{An \(L^1\) bound for Ricci curvature on K\"ahler manifolds}
\author[W. Jiang]{Wenshuai Jiang}
\address{School of Mathematical Sciences, Zhejiang University, Hangzhou 310058, China}
\email{wsjiang@zju.edu.cn}
\author[G. Liu]{Gang Liu}
\address{School of Mathematical Sciences, East China Normal University, Shanghai 200241, China}
\email{gliu@math.ecnu.edu.cn}
\date{October 20, 2026}

\begin{document}

\begin{abstract}
We prove a uniform \(L^1\) bound for Ricci curvature on
noncollapsed K\"ahler manifolds with Ricci curvature bounded below.
The proof combines the codimension two neck decomposition and entropy
estimates with local K\"ahler geometry on regular balls. The key step is
a linear control of the Ricci integral on a regular ball by the entropy
pinching, obtained from a decomposition of the scalar-curvature integral
that exploits the K\"ahler structure.
\end{abstract}
\maketitle

\setcounter{tocdepth}{2}
\tableofcontents
\clearpage

\section{Introduction}
\label{sec:introduction}

In this paper, we prove a uniform \(L^1\) bound for Ricci curvature
on noncollapsed K\"ahler manifolds with Ricci curvature bounded
below. Our result establishes the K\"ahler case of the
noncollapsing form of Yau's scalar-curvature conjecture
\cite{Yau1992}, as stated in
\cite[Conjecture~2.18]{Naber2020}.

\begin{theorem}\label{thm:main}
Let \((M^n,g,J,p)\) be a complete K\"ahler manifold of real
dimension \(n=2m\ge2\), with \(\Ric_g\ge-(n-1)g\) and
\(\Vol_g(B_1(p))\ge v>0\); then
\begin{equation}\label{eq:main}
 \int_{B_1(p)}|\Ric_g|\dd V_g\le C(n,v).
\end{equation}
\end{theorem}

Under a sectional-curvature lower bound, a local curvature integral
estimate was established by Petrunin \cite{Petrunin2009}, with an
alternative proof given by Li \cite{Li2025}. The problem is more
difficult under only a Ricci lower bound. The counterexamples of Hao and Zhu
\cite{HaoZhu2026} and Xu \cite{Xu2026Collapse} disprove the
asymptotic form of Yau's conjecture, while Xu's collapsing examples
also show that a uniform local bound cannot hold on general
Riemannian manifolds without a noncollapsing assumption.
Under a noncollapsing assumption, Jiang and Naber obtained \(L^p\)
estimates for \(0<p<1\) in the unpublished work
\cite{JiangNaberUnpublished}; see also Naber~\cite{Naber2020}.
A proof was later given by Cucinotta and Mondino
\cite{CucinottaMondino2026}. Here we study the case \(p=1\).
Other related results on scalar-curvature
integrals and Cohn--Vossen-type inequalities can be found in
\cite{ChenXuZhang2026,Deng2026,Ma2026,MunteanuWang2025,Xu2020,Xu2024,Yang2013,Zhu2022Comparison,Zhu2022}.

K\"ahler geometry provides additional structure for proving the
\(L^1\) bound. Under positive or nonnegative bisectional curvature,
the work of Chen and Zhu \cite{ChenZhu2005}, Ni \cite{Ni2005},
and Liu \cite{Liu2024} relates volume growth to scalar-curvature
decay and integral estimates. Related developments
include Liu's Cohn--Vossen-type inequality for certain noncompact
K\"ahler manifolds \cite{LiuCohnVossen2022} and the Ricci-form
estimates of Liu and Ni \cite{LiuNi2026} under additional
nonnegative-curvature and Ricci-pinching assumptions.

For the local problem under a Ricci lower bound, the results of
Liu and Sz\'ekelyhidi \cite{LiuSzekelyhidi2022} are particularly
relevant. On almost Euclidean balls, they constructed holomorphic
coordinates and obtained local scalar-curvature integral estimates
without any polarization assumption. For general noncollapsed
unit balls, they proved the corresponding integral bound in the
polarized setting. Under a Ricci lower bound, these scalar-curvature
integral bounds are equivalent to local \(L^1\) bounds for Ricci
curvature. Theorem~\ref{thm:main} removes the polarization
assumption from the estimate on general balls.

Our proof of this extension combines the almost Euclidean
estimates with an analysis of regions close to singular cones.
The foundational work of Cheeger and Colding
\cite{CheegerColding1996,CheegerColding1997,CheegerColding2000II,CheegerColding2000III}
established the almost splitting theorem and the structure theory
of noncollapsed limit spaces, providing the basic analytic and
geometric framework for this approach. In the K\"ahler setting,
the fundamental work of Cheeger, Colding, and Tian
\cite{CheegerColdingTian2002} showed how the complex structure
gives stronger control of the singularities of limit spaces.
In particular, they used the parallel complex structure to obtain
almost holomorphic functions from real splitting functions.
This idea also appears in \cite{LiuSzekelyhidi2022} and is used
in the complex pairing argument in Lemma~\ref{lem:complex-pairing}.
Quantitative control of the singular geometry was developed by
Cheeger and Naber \cite{CheegerNaber2013}, who introduced
quantitative stratification and obtained volume estimates for
quantitative singular strata. For noncollapsed manifolds with
bounded Ricci curvature, they also proved codimension four
regularity and \(L^q\) curvature estimates for \(q<2\)
\cite{CheegerNaber2015}.

Building on this work, Jiang and Naber~\cite{JiangNaber2021}
proved the \(L^2\) bound for the full curvature tensor under a two-sided
Ricci bound and a noncollapsing assumption. Their proof decomposes
a ball into regular regions and codimension four necks, with
uniform packing estimates. On each neck, estimates for splitting
and cone functions give a local \(L^2\) curvature estimate that is linear
in the drop of a monotone volume quantity. Monotonicity then allows
these estimates to be summed over all scales. Combining the neck
estimates with the bounds on the regular regions then gives
the uniform \(L^2\) curvature bound.
For applications of this framework and related work on Yang--Mills
fields, harmonic maps, geometric flows and parabolic equations, see
\cite{HuangZhang2026,HuangJiang2024,HuangJiang2025,NaberValtorta2019,NaberValtorta2026}.

Our proof follows this strategy in the codimension two setting,
using the neck structure, decomposition, and entropy estimates
of Cheeger, Jiang, and Naber~\cite{CheegerJiangNaber2021}.
The local analysis is more difficult here: with only a Ricci lower
bound, we no longer have the uniform regularity in harmonic coordinates
used in the \(L^2\) argument of~\cite{JiangNaber2021}.
The almost Euclidean balls in this decomposition are controlled
by the local estimates of Liu and Sz\'ekelyhidi~\cite{LiuSzekelyhidi2022}.
The remaining difficulty is to control the Ricci integral on
regular balls inside a neck. To sum these local bounds over all
scales, we need an estimate that is linear in the entropy pinching.

This linear bound is the new estimate in
Proposition~\ref{prop:regular}. Using a section of the canonical
bundle and the cone function, we split the scalar-curvature
integral into two terms. We control the weighted term by Ricci
contractions through Bochner formulas. For the other term, we use
products of two \(L^2\) errors of square-root size and squared
derivative terms; a holomorphic correction and integration by parts
then give the required linear bound.
Proposition~\ref{prop:neck-energy} converts this estimate into
linear entropy control on regular balls, which can be summed
using entropy monotonicity and the neck packing estimates.

Section~\ref{sec:neck-background} records the needed results and
consequences from \cite{CheegerJiangNaber2021,LiuSzekelyhidi2022}.
Section~\ref{sec:regular} proves the key local linear estimate; the
development of this argument was substantially assisted by AI, as
explained in the Acknowledgments. Section~\ref{sec:proof} follows the
Jiang--Naber neck summation scheme to complete the proof.

\section{Neck regions and entropy estimates}
\label{sec:neck-background}

We recall the geometric decomposition and localized entropy estimates from
\cite{CheegerJiangNaber2021}, together with the local K\"ahler results from
\cite{LiuSzekelyhidi2022} used in the proof. We also record several
consequences of these results in the form needed below.
Theorem~\ref{thm:neck-functions} provides a splitting map and a cone
function whose Hessian energies and Ricci contractions are controlled
linearly by entropy pinching;
Proposition~\ref{prop:neck-energy} collects the estimates needed for
the curvature calculation.
All manifolds below are complete and smooth.

\subsection{Symmetry, volume ratios, and neck regions}

We first describe the cone models and centers of a codimension two
neck. The Euclidean directions in the cone model determine the
axis in the neck definition.

\begin{definition}
\label{def:symmetry}
A ball \(B_s(x)\) is \((k,\eps)\)-symmetric if there is a metric cone
\[
 X'=\mathbb R^k\times C(Z),\qquad x_0=(0,o),
\]
where \(o\) is the cone vertex, such that
\[
 d_{\GH}\bigl(B_s(x),B_s(x_0)\bigr)<\eps s.
\]
Given a corresponding GH approximation
\(\iota:B_s(x_0)\longrightarrow B_s(x)\), write
\[
 L_{x,s}=\iota\bigl((\mathbb R^k\times\{o\})\cap B_s(x_0)\bigr).
\]
We then say that the ball is symmetric with respect to \(L_{x,s}\).
The set \(L_{x,s}\) is the image of the model axis under this
approximation.
\end{definition}

For \(\Ric\ge-(n-1)\kappa g\), \(\kappa\ge0\), write
\begin{equation}\label{eq:volume-ratio}
 \V_s^\kappa(x)=\frac{\Vol(B_s(x))}{\Vol_{-\kappa}(B_s)},
\end{equation}
where the denominator is the volume in the simply connected
space form of curvature \(-\kappa\). We omit \(\kappa\) when the
lower bound is fixed.

A neck is the region outside a family of center balls, where the
space is close to codimension two cones at each scale. The conditions
below are used in the packing and function estimates.

\begin{definition}[{\cite{CheegerJiangNaber2021}}]
\label{def:neck}
Let \((M^n,g)\) be a complete smooth Riemannian manifold,
with \(R,v,\delta_0,\eta>0\) and \(\kappa\ge0\).
Suppose \(\Ric\ge-(n-1)\kappa g\),
\(\kappa R^2\le\delta_0^2\), and
\(\Vol(B_R(q))\ge vR^n\). Let \(\C\subset B_{2R}(q)\) be closed,
with \(q\in\C\), and let \(x\mapsto r_x>0\) be continuous on \(\C\).
The region
\begin{equation}\label{eq:neck}
 \N=B_{2R}(q)\setminus\overline B_{r_x}(\C),\qquad
 \overline B_{r_x}(\C)=\bigcup_{x\in\C}\overline B_{r_x}(x),
\end{equation}
is an \((n-2,\delta_0,\eta)\)-neck if:
\begin{enumerate}
\item[\textup{(n1)}]
The balls \(B_{\tau_n^2r_x}(x)\), \(x\in\C\), are pairwise
 disjoint and contained in \(B_{2R}(q)\), where
\(\tau_n=10^{-10n}\omega_n\) and \(\omega_n=\Vol(B_1(0^n))\).
\item[\textup{(n2)}]
For every \(x\in\C\),
\[
 \left|\V_{\delta_0^{-1}R}(x)-\V_{\delta_0r_x}(x)\right|
       <\delta_0^2.
\]
\item[\textup{(n3)}]
For \(x\in\C\) and \(r_x\le s\le\delta_0^{-1}R\),
\(B_s(x)\) is \((n-2,\delta_0^2)\)-symmetric with respect to
\(L_{x,s}\), but is not \((n-1,\eta)\)-symmetric.
\item[\textup{(n4)}]
For \(x\in\C\), \(s\ge r_x\), and
\(B_{2s}(x)\subset B_{2R}(q)\),
\[
 L_{x,s}\subset B_{\tau_ns}(\C),\qquad
 \C\cap B_s(x)\subset B_{\tau_ns}(L_{x,s}).
\]
\item[\textup{(n5)}]
\(\Lip r_x\le\delta_0\).
\end{enumerate}
Here \(B_t(D)=\bigcup_{z\in D}B_t(z)\).
\end{definition}

All center radii are positive in the smooth case.
We weight each center by the \((n-2)\)-power of its radius;
this is the measure used to average entropy pinching and to sum
the local curvature estimates.

\begin{definition}
\label{def:packing}
Let \(\N=B_{2R}(q)\setminus\overline B_{r_x}(\C)\) be a
smooth \((n-2,\delta_0,\eta)\)-neck. Its packing measure is
\begin{equation}\label{eq:mu}
 \mu=\mu_{\N}=\sum_{x\in\C}r_x^{n-2}\delta_x,
\end{equation}
where \(\delta_x\) is the unit Dirac mass at \(x\).
\end{definition}

\subsection{Neck structure and neck decomposition}

The following two results control the size of individual necks and
the total size of the decomposition. The Ahlfors estimate for the
center measure will be used to sum the local entropy bounds on
each neck.

\begin{theorem}[Neck Structure Theorem {\cite{CheegerJiangNaber2021}}]
\label{thm:neck-structure}
Let \(\N=B_{2R}(q)\setminus\overline B_{r_x}(\C)\) be a
smooth \((n-2,\delta_0,\eta)\)-neck in a complete
Riemannian manifold, with \(\Vol(B_R(q))\ge vR^n\).
If \(\delta_0\le\delta_0(n,v,\eta)\), its packing measure
\(\mu=\sum_{x\in\C}r_x^{n-2}\delta_x\) satisfies
\begin{equation}\label{eq:ahlfors}
 A_0^{-1}s^{n-2}\le\mu(B_s(x))\le A_0s^{n-2},
 \qquad A_0=A_0(n),
\end{equation}
whenever \(x\in\C\), \(s\ge r_x\), and
\(B_{2s}(x)\subset B_{2R}(q)\).
\end{theorem}

The next theorem provides the covering used in the final step of
the proof. The bound on the sum of the \((n-2)\)-powers of the
radii matches the scaling of both local Ricci estimates.

\begin{theorem}[Neck Decomposition Theorem {\cite{CheegerJiangNaber2021}}]
\label{thm:neck-decomposition}
Let \((M^n,g,p)\) be a complete smooth Riemannian manifold with
\(\Ric\ge-(n-1)g\) and
\(\Vol(B_1(p))\ge v>0\). For \(\eta>0\) and
\(\delta_0\le\delta_0(n,v,\eta)\), there is a countable covering
\begin{equation}\label{eq:neck-decomposition}
 B_1(p)\subset
 \bigcup_a\bigl(\N_a\cap B_{r_a}(x_a)\bigr)
       \ \cup\ \bigcup_b B_{r_b}(x_b),
\end{equation}
where \(0<r_a,r_b\le1\), each
\(\N_a\subset B_{2r_a}(x_a)\) is an
\((n-2,\delta_0,\eta)\)-neck, each \(B_{2r_b}(x_b)\) is
\((n-1,2\eta)\)-symmetric, and
\begin{equation}\label{eq:packing-content}
 \sum_a r_a^{n-2}+\sum_b r_b^{n-2}
       \le C(n,v,\delta_0,\eta).
\end{equation}
\end{theorem}

We use codimension one rigidity to regard the latter balls as
almost Euclidean on a fixed interior scale.
The codimension two cone models are
\(\mathbb R^{n-2}\times C(S^1_\alpha)\), with
\(\alpha\ge\alpha(n,v)>0\); they are locally Euclidean away
from their axes.
These give the usual regular-ball geometry of the necks.

\subsection{The localized pointed entropy}
\label{sec:entropy}

We use the localized pointed entropy to measure the error in the
cone approximation at each scale. Its monotonicity and bounded
range will control the sum of entropy differences in the neck estimate.

For this subsection assume
\begin{equation}\label{eq:entropy-geometry}
 \Ric_g\ge-(n-1)\delta_1^2g,\qquad
 \Vol(B_s(q))\ge vs^n\quad(0<s\le\delta_1^{-1}).
\end{equation}

We use the same cutoff for all centers. The correction terms make
the entropy monotone under the Ricci lower bound.

\begin{definition}[{\cite{CheegerJiangNaber2021}}]
\label{def:entropy}
Let \((M^n,g,q)\) be a complete smooth Riemannian manifold,
and let \(v,\delta_1>0\). Assume
\[
 \Ric_g\ge-(n-1)\delta_1^2g,\qquad
 \Vol(B_\rho(q))\ge v\rho^n\quad(0<\rho\le\delta_1^{-1}).
\]
Write the heat kernel and its measure as
\[
 \rho_t(x,y)=(4\pi t)^{-n/2}e^{-f_{x,t}(y)},\qquad
 \rho_t(x,dy)=\rho_t(x,y)\dd V_g(y).
\]
Choose \(0\le\phi\le1\), with
\[
 \phi=1\ \text{on }B_{\delta_1^{-1}}(q),\qquad
 \operatorname{supp}\phi\subset B_{2\delta_1^{-1}}(q),\qquad
 |\Delta\phi|+|\nabla\phi|^2\le C(n)\delta_1^2.
\]
For fixed \(x,t\), put
\[
 W_t=2t\Delta_f f+t|\nabla f|^2+f-n,
 \qquad
 \Delta_f=\Delta-\langle\nabla f,\nabla\cdot\rangle,
 \qquad f=f_{x,t}.
\]
For \(x\in B_{\delta_1^{-1}/2}(q)\) and
\(0<t\le\delta_1^{-2}\), define
\begin{equation}\label{eq:Wdefinition}
\begin{split}
 \W_t^{\delta_1}(x)
 ={}&\int_M W_t\phi\,\rho_t(x,dy)\\
 &-4(n-1)\delta_1^2\int_0^t
             s\int_M|\nabla f_{x,s}|^2\phi\,\rho_s(x,dy)\dd s\\
 &-C(n,v)\delta_1^2\int_0^t
             \exp\!\left(-\frac{1}{100\delta_1^2s}\right)\dd s.
\end{split}
\end{equation}
The constant \(C(n,v)\) is chosen to ensure monotonicity. A single cutoff and the same constant are retained
whenever entropy values are compared.
\end{definition}

The next theorem relates the entropy drop to squared Hessian
errors and Ricci contractions. It also compares entropy with the
logarithmic volume ratio when the latter changes little across scales.

\begin{theorem}[{\cite{CheegerJiangNaber2021}}]
\label{thm:entropy}
Let \((M^n,g,q)\) be a complete smooth Riemannian manifold
and suppose, for \(v,\delta_1>0\), that
\[
 \Ric_g\ge-(n-1)\delta_1^2g,\qquad
 \Vol(B_\rho(q))\ge v\rho^n\quad(0<\rho\le\delta_1^{-1}).
\]
Use the localized entropy \(\W^{\delta_1}\) with one fixed
cutoff based at \(q\). For
\(x\in B_{\delta_1^{-1}/2}(q)\) and
\(0<t\le\delta_1^{-2}\), one has
\begin{equation}\label{eq:monotone}
\begin{split}
 \frac{d}{dt}\W_t^{\delta_1}(x)
 \le{}&-2t\int_M\bigg(
      \left|\Hess f_{x,t}-\frac{g}{2t}\right|^2
      +\Ric(\nabla f_{x,t},\nabla f_{x,t})\\
 &\hspace{39mm}+2(n-1)\delta_1^2|\nabla f_{x,t}|^2
                    \bigg)\phi\,\rho_t(x,dy)\le0.
\end{split}
\end{equation}
Moreover, given \(\eps>0\), if
\(\delta_1\le\delta_1(n,v,\eps)\), \(0<t\le10\), and
\begin{equation}\label{eq:entropy-volume-pinching}
 \left|\V^{\delta_1^2}_{\delta_1^{-1}\sqrt t}(x)
       -\V^{\delta_1^2}_{\delta_1\sqrt t}(x)\right|\le\delta_1,
\end{equation}
then
\begin{equation}\label{eq:entropy-volume}
 \left|\W_t^{\delta_1}(x)
           -\log\V^{\delta_1^2}_{\sqrt t}(x)\right|\le\eps.
\end{equation}
\end{theorem}

The bounds below and monotonicity allow us to sum the entropy
drops over all scales.

\begin{remark}\label{rem:entropy-bound}
The heat kernel estimates and monotonicity give, for the same centers and times,
\begin{equation}\label{eq:Wrange}
 -C(n,v)\le\W_t^{\delta_1}(x)\le0,
 \qquad \lim_{t\downarrow0}\W_t^{\delta_1}(x)=0.
\end{equation}
\end{remark}

Under \(\widetilde g=a^{-2}g\), keeping the cutoff function \(\phi\) fixed,
\begin{equation}\label{eq:Wscale}
 \W_{a^2t}^{\delta_1,g,\phi}(x)
       =\W_t^{a\delta_1,\widetilde g,\phi}(x).
\end{equation}

\subsection{The splitting map and the cone function}

We next recall the functions used in the curvature calculation:
a harmonic splitting map and a Poisson approximation to squared
distance. Entropy pinching at independent centers controls the
errors in all the splitting directions.

\begin{definition}
\label{def:pinching}
For \(k\ge1\), a collection
\(U=\{x_0,\ldots,x_k\}\subset B_{2r}(x)\) is
\((k,\alpha)\)-independent at scale \(r\) if
\[
 d_{\GH}(U,U')\ge\alpha r
\]
for every configuration
\(U'=\{x'_0,\ldots,x'_k\}\subset\mathbb R^{k-1}\), with repeated
points allowed. The GH distance is taken between the resulting
finite sets with their induced metrics.
For \(2r^2\le\delta_1^{-2}\), the
\((k,\alpha,\delta_1)\)-entropy pinching is
\begin{equation}\label{eq:k-pinching}
 E_r^{k,\alpha,\delta_1}(x)
  =\inf_{\{x_0,\ldots,x_k\}}
       \sum_{a=0}^k
       \left(\W_{r^2}^{\delta_1}(x_a)
                         -\W_{2r^2}^{\delta_1}(x_a)\right),
\end{equation}
where the infimum is over these independent collections in
\(B_{2r}(x)\), with each center in
\(B_{\delta_1^{-1}/2}(q)\). All entropy values use one fixed
localization based at \(q\).
\end{definition}

The harmonic splitting map represents the Euclidean factor of
the cone model. Its defining estimates express approximate
orthonormality and small Hessian energy.

\begin{definition}
\label{def:splitting}
A map \(u=(u^1,\ldots,u^k):B_r(x)\to\mathbb R^k\) is a harmonic
\((k,\eps)\)-splitting map if \(\Delta u^a=0\) for every \(a\), and
\begin{gather*}
 \max_{1\le a\le k}\sup_{B_r(x)}|\nabla u^a|\le1+\eps,\\
 \fintt_{B_r(x)}
     |\langle\nabla u^a,\nabla u^b\rangle-\delta_{ab}|\dd V
       \le\eps\qquad(1\le a,b\le k),\\
 r^2\fintt_{B_r(x)}\sum_{a=1}^k|\Hess u^a|^2\dd V
       \le\eps^2.
\end{gather*}
\end{definition}

We next pass from pinching averaged over the neck centers to
pinching at an independent collection of centers. The same
center-selection argument appears in Cheeger--Jiang--Naber and
Jiang--Naber. We record the resulting consequence and omit the proof.

\begin{remark*}
Consider a smooth \((n-2,\delta_0,\eta)\)-neck of dimension
\(n\ge3\) satisfying
\eqref{eq:entropy-geometry}, with
\(\delta_1R\le4\delta_0\le1/100\). Let
\[
 y\in\N\cap B_R(q),\qquad
 r=d(y,\C)=d(y,z)>0,\qquad
 z\in\C,\qquad s=Lr\le10^{-6}R,
\]
and put \(k=n-2\) and
\begin{equation}\label{eq:average-pinching}
 P=s^{-k}\int_{\C\cap B_s(z)}
       \left(\W_{s^2/2}^{\delta_1}(x)
                   -\W_{2s^2}^{\delta_1}(x)\right)\dd\mu(x).
\end{equation}
For \(L\ge L_0(n,v,\eta)\ge32\) and
\(0<\delta_0\le\delta_0(n,v,\eta,L)\), the Ahlfors estimate
and the center-selection argument give \((k,\alpha)\)-independent centers
\(x_0,\ldots,x_k\in\C\cap B_s(z)\), at scale \(s\), such that
\begin{equation}\label{eq:selected-centers}
 \sum_{j=0}^k
       \left(\W_{s^2/2}^{\delta_1}(x_j)
                    -\W_{2s^2}^{\delta_1}(x_j)\right)\le CP,
 \qquad \alpha=\alpha(n)>0.
\end{equation}
Here \(C=C(n,v,\eta,L)\); all entropy values retain the
cutoff and parameter fixed for the neck.
\end{remark*}

The next theorem provides a splitting map and a cone function
with linear bounds for their Hessian energies and Ricci contractions.
After rescaling, the following result follows from the sharp
cone-splitting and Poisson regularization estimates.

\begin{theorem}[{\cite{CheegerJiangNaber2021}}]
\label{thm:neck-functions}
Given \(n\ge3\) and \(v,\eta,\eps,\alpha>0\), let
\(\delta_0\le\delta_0(n,v,\alpha,\eps)\), and let
\[
 \N=B_{2R}(q)\setminus\overline B_{r_x}(\C)
\]
be a smooth \((n-2,\delta_0,\eta)\)-neck with
\(r_q\le R\), in a complete Riemannian manifold \((M^n,g)\).
Set \(\delta_1=4\delta_0/R\), and assume
\[
 \Ric_g\ge-(n-1)\delta_1^2g,\qquad
 \Vol(B_\rho(q))\ge v\rho^n\quad(0<\rho\le\delta_1^{-1}).
\]
Use the localized entropy with one fixed cutoff based at \(q\).
There exist a harmonic
\((n-2,\eps)\)-splitting map
\(u=(u^1,\ldots,u^{n-2}):B_{2R}(q)\to\mathbb R^{n-2}\)
and a function \(h:B_{2R}(q)\to\mathbb R\) with the following
properties:
\begin{equation}\label{eq:splitting-energy}
\begin{split}
 R^2\fintt_{B_{2R}(q)}\sum_{a=1}^{n-2}\bigl(
       |\Hess u^a|^2+\Ric(\nabla u^a,\nabla u^a)
       +2(n-1)\delta_1^2|\nabla u^a|^2\bigr)\dd V\\
 \le C(n,v,\alpha)E_R^{n-2,\alpha,\delta_1}(q),
\end{split}
\end{equation}
\begin{equation}\label{eq:poisson-equation}
 \Delta h=2n\qquad\text{on }B_{2R}(q),
\end{equation}
\begin{equation}\label{eq:poisson-linear}
\begin{split}
 \fintt_{B_{2R}(q)}\bigl(
       |\Hess h-2g|^2+\Ric(\nabla h,\nabla h)
       +2(n-1)\delta_1^2|\nabla h|^2\bigr)\dd V\\
 \le C(n,v)
       \left(\W_{R^2}^{\delta_1}(q)-\W_{2R^2}^{\delta_1}(q)\right),
\end{split}
\end{equation}
\begin{equation}\label{eq:poisson-eikonal}
 \fintt_{B_{2R}(q)}\bigl||\nabla h|^2-4h\bigr|^2\dd V
 \le C(n,v)R^4
       \left(\W_{R^2}^{\delta_1}(q)-\W_{2R^2}^{\delta_1}(q)\right),
\end{equation}
and
\begin{equation}\label{eq:poisson-sup}
 \sup_{B_{2R}(q)}|\nabla h|\le C(n,v)R,\qquad
 \sup_{B_{2R}(q)}|h-d_q^2|\le\eps R^2.
\end{equation}
\end{theorem}

The neck hypotheses provide the noncollapsing needed for this
localization, after decreasing \(v\) by a dimensional factor.
For the local curvature calculation, we need one choice of
\(u,h\) whose errors are controlled by the averaged pinching
\(P\), with \(h-|u|^2\) positive on a smaller ball centered at \(y\).

\begin{corollary}\label{cor:common-functions}
Let \(\N\subset B_{2R}(q)\) be a smooth
\((n-2,\delta_0,\eta)\)-neck in a complete Riemannian
manifold of dimension \(n\ge3\). Suppose, for \(v,\delta_1>0\),
that
\[
 \Ric_g\ge-(n-1)\delta_1^2g,\qquad
 \Vol(B_\rho(q))\ge v\rho^n\quad(0<\rho\le\delta_1^{-1}),
 \qquad \delta_1R\le4\delta_0\le1/100.
\]
Fix \(L\ge L_0(n,v,\eta)\ge32\), and take
\[
 y\in\N\cap B_R(q),\qquad z\in\C,\qquad
 r=d(y,\C)=d(y,z)>0,\qquad s=Lr\le10^{-6}R.
\]
Let \(\mu=\mu_{\N}\) be the packing measure, and define
\(P\) by \eqref{eq:average-pinching}, using one fixed entropy
cutoff based at \(q\). Put
\[
 k=n-2,\qquad \kappa=(n-1)\delta_1^2,\qquad
 \mathcal R=\Ric+\kappa g\ge0.
\]
There are \(\eps_0>0\), \(C<\infty\), and
\(0<c\le C<\infty\), all depending only on
\(n,v,\eta,L\), such that, if \(P+\delta_0\le\eps_0\),
there are harmonic functions \(u^1,\ldots,u^k\) and a function \(h\)
on \(B_{8s}(z)\), with
\[
 \Delta h=2n,\qquad
 |u|\le Cs,\quad |du|\le C,\quad |h|\le Cs^2,\quad |dh|\le Cs,
\]
such that
\begin{equation}\label{eq:common-energy}
\begin{split}
 s^{-n}\int_{B_{8s}(z)}\Bigl[
 &\sum_{a,b=1}^k
     |\langle du^a,du^b\rangle-\delta_{ab}|^2\\
 &+s^2\sum_{a=1}^k
       \bigl(|\Hess u^a|^2+\mathcal R(\nabla u^a,\nabla u^a)\bigr)\\
 &+|\Hess h-2g|^2+\mathcal R(\nabla h,\nabla h)
           +s^{-4}(|dh|^2-4h)^2\Bigr]\dd V_g\le CP.
\end{split}
\end{equation}
After adding constants to \(u^a\), we also have
\begin{equation}\label{eq:homogeneity}
 \ell_a=\langle dh,du^a\rangle-2u^a,\qquad
 s^{-n-2}\int_{B_{2s}(z)}\sum_{a=1}^k\ell_a^2\dd V_g\le CP.
\end{equation}
For \(b=h-|u|^2\),
\begin{equation}\label{eq:transverse-positive}
 cs^2\le b\le Cs^2\qquad\text{on }B_{r/8}(y).
\end{equation}
\end{corollary}

\begin{proof}
Use the fixed-cutoff local construction underlying
Theorem~\ref{thm:neck-functions}. Apply it to the selected centers on
a common fixed enlargement of \(B_{8s}(z)\), choosing \(h\)
with radial center \(x_0\). The cutoff, entropy parameter,
and time interval \([s^2/2,2s^2]\) remain fixed. Constant
Gram normalization and the Poincar\'e inequality give
\eqref{eq:common-energy}; the enlargement changes only the constants.

To prove \eqref{eq:homogeneity}, differentiate:
\[
 d\ell_a=(\Hess h-2g)(\nabla u^a,\cdot)
                   +\Hess u^a(\nabla h,\cdot).
\]
Replace \(u^a\) by
\(u^a+\tfrac12\fintt_{B_{2s}(z)}\ell_a\dd V_g\).
The new \(\ell_a\) has mean zero, so
\[
\begin{split}
 \int_{B_{2s}(z)}\sum_a\ell_a^2\dd V_g
 &\le Cs^2\int_{B_{4s}(z)}\sum_a|d\ell_a|^2\dd V_g\\
 &\le Cs^2\int_{B_{4s}(z)}
       \left(|\Hess h-2g|^2+s^2\sum_a|\Hess u^a|^2\right)\dd V_g
 \le Cs^{n+2}P.
\end{split}
\]
In the normalized cone limit, \(h=|Z-Z_0|^2+\rho^2\) and
\(u=OZ+c\); the mean-zero condition gives \(c=-OZ_0\).
Thus \(b\to\rho^2\) uniformly. Since \(\rho\ge r/4\)
on \(B_{r/8}(y)\), the conic approximation gives
\eqref{eq:transverse-positive} after \(L\) is fixed.
\end{proof}

\subsection{Holomorphic charts and Ricci estimates on almost Euclidean balls}
\label{sec:kahler-background}

The preceding results apply to general Riemannian manifolds.
We now state the holomorphic chart and Ricci integral estimates
used below.
In this subsection, let \((M^n,g,J)\) be a K\"ahler manifold
of real dimension \(n=2m\ge4\), with complex structure \(J\) and
K\"ahler form \(\omega(X,Y)=g(JX,Y)\). Thus \(\nabla J=0\).
Write \(S=\operatorname{tr}_{\mathbb R}\Ric_g\). In holomorphic
coordinates \(z=(z_1,\ldots,z_m)\), we use
\[
 \omega=\frac{i}{2}\sum_{a,b}g_{a\bar b}\,dz_a\wedge d\bar z_b,
 \qquad
 \Ric(\omega)=-\frac{i}{2\pi}\partial\bar\partial
                                      \log\det(g_{a\bar b}).
\]
With this normalization,
\[
 2\pi\Ric(\omega)(X,Y)=\Ric_g(JX,Y).
\]
The symbol \(\Ric_g\), or \(\Ric\) without the argument \(\omega\),
always denotes the real Ricci tensor. We denote the Levi--Civita
connection by \(\nabla\) and its induced Chern covariant derivative
on holomorphic tensor bundles by \(D=D'+D''\), where
\(D''=\bar\partial\). Here \(D'\partial f=(\Hess f)^{2,0}\)
for a smooth function \(f\).

We also write
\[
 \omega_0=\frac{i}{2}\sum_jdz_j\wedge d\bar z_j,
 \qquad dV_0=\frac{\omega_0^m}{m!}.
\]
Here \(B_\rho(y)\) denotes a metric ball, and \(B^0_\rho(0)\)
a Euclidean coordinate ball, identified with its inverse image
in the chart. Coordinate norms have subscript \(0\); all other
norms use \(g\). We keep these conventions in
Section~\ref{sec:regular}.

The next theorem gives holomorphic coordinates and a bounded
K\"ahler potential on almost Euclidean balls. We use these in
Section~\ref{sec:regular}.

\begin{theorem}[{\cite{LiuSzekelyhidi2022}}]
\label{thm:holomorphic-chart}
For \(n=2m\ge4\), there is \(\eps_0(n)>0\) with the following
property. Let \(0<\eps\le\eps_0(n)\), \(r>0\), and suppose that
\(B_{\eps^{-1}r}(x)\) is relatively compact in a smooth
K\"ahler manifold of real dimension \(n=2m\ge4\).
Assume, on this ball, that
\[
 \Ric_g>-\eps r^{-2}g,\qquad
 d_{\GH}\bigl(B_{\eps^{-1}r}(x),
                    B_{\eps^{-1}r}(0^n)\bigr)<\eps r.
\]
Then there is a holomorphic coordinate map
\(F:B_r(x)\longrightarrow\mathbb C^m\) that is a
\(\Psi(\eps\mid n)r\)-GH approximation to its image, where
\(\Psi(\eps\mid n)\to0\) as \(\eps\downarrow0\).
There is also a smooth real potential \(\rho\) on \(B_r(x)\)
such that
\[
 \omega=i\partial\bar\partial\rho,
 \qquad \sup_{B_r(x)}|\rho|\le C(n)r^2.
\]
The last estimate is the bounded-potential consequence of the
potential approximation in the cited theorem.
\end{theorem}

The scalar-curvature estimate of Liu and Sz\'ekelyhidi also gives
the following Ricci \(L^1\) estimate. This follows from
\(|\Ric|\le S+2n\kappa\) when \(\Ric\ge-\kappa g\),
together with volume comparison.

\begin{theorem}[{\cite{LiuSzekelyhidi2022}}]
\label{thm:scalar-almost-euclidean}
Given \(n=2m\ge4\) and \(v>0\), there are
\(\eps=\eps(n,v)>0\) and \(C=C(n,v)<\infty\)
with the following property. Let \((M^n,g,J)\) be a complete
K\"ahler manifold of real dimension \(n\), and let \(r>0\).
If
\[
\begin{gathered}
 \Ric_g\ge-(n-1)\eps r^{-2}g,\qquad
 \Vol_g(B_r(x))\ge vr^n,\\
 d_{\GH}\bigl(B_{4r}(x),B_{4r}(0^n)\bigr)<\eps r,
\end{gathered}
\]
then
\[
 \int_{B_r(x)}|\Ric_g|\dd V_g\le Cr^{n-2}.
\]
\end{theorem}

Since the neck models are flat away from their axes,
Theorem~\ref{thm:holomorphic-chart} and the interior gradient
estimate for its harmonic coordinate functions apply on regular
balls. The following corollary gives the chart, potential, and
metric bounds needed for the holomorphic correction and curvature
integral in Section~\ref{sec:regular}.

\begin{corollary}\label{cor:neck-chart}
For \(n=2m\ge4\), \(v,\eta>0\), and \(0<\eps<1/100\),
there are \(a_0\in(0,1/8)\), \(\delta>0\), and
\(C<\infty\), depending only on \(n,v,\eta,\eps\), with the
following property.
Let \(\N\subset B_{2R}(q)\) be a smooth
\((n-2,\delta_0,\eta)\)-neck
in a complete K\"ahler manifold of real dimension \(n=2m\ge4\),
with \(\delta_0\le\delta\) and
\(\Vol_g(B_R(q))\ge vR^n\).
For \(y\in\N\cap B_R(q)\), put \(r=d(y,\C)>0\).
There is a holomorphic coordinate map
\[
 z=(z_1,\ldots,z_m):B_{a_0r}(y)\longrightarrow\mathbb C^m,
 \qquad z(y)=0,
\]
which, for the rescaled metric \(\widetilde g=r^{-2}g\),
is an \(\eps a_0\)-GH approximation to its image.
In these coordinates \(\widetilde g\) admits a K\"ahler
potential \(\widetilde\rho\) and satisfies
\[
 |\widetilde\rho|\le C,\qquad
 \widetilde g\ge C^{-1}z^*g_{\mathbb C^m}.
\]
\end{corollary}

For the application in Section~\ref{sec:proof}, choose the chart accuracy once
in terms of \(n,v,\eta\), small enough for
Theorem~\ref{thm:scalar-almost-euclidean} and the interior ball
inclusions below. Decrease the associated radius constant by
a factor of \(8\), and retain the notation \(a_0\).
This fixes \(a_0=a_0(n,v,\eta)<1/64\) and a chart on
\(B_{8a_0r}(y)\) whose distance distortion for \(r^{-2}g\)
is at most \(10^{-3}a_0\).
The resulting constants also depend only on
\(n,v,\eta\).

\subsection{Complex pairing of the splitting functions}

We combine pairs of real splitting functions into nearly holomorphic
functions and collect their estimates at the scale of a regular ball.
Proposition~\ref{prop:neck-energy} provides the energy bounds
used to apply the local estimate in Section~\ref{sec:proof}.

The next lemma gives a linear bound for the squared
antiholomorphic derivatives after a constant orthogonal change
of the real splitting functions. We state it at unit scale;
the following proposition applies it at scale \(s=Lr\).

\begin{lemma}\label{lem:complex-pairing}
Fix \(n=2m\ge4\), \(v,\eta>0\), and \(C_0\ge1\), and put
\(k=n-2\). Let
\[
 \N=B_{2R}(q)\setminus\overline B_{r_x}(\C)
\]
be an \((n-2,\delta_0,\eta)\)-neck in a complete
K\"ahler manifold \((M^n,g,J)\) of real dimension \(n\).
Suppose that \(z\in\C\) and
\[
 r_z\le1/32,\qquad B_{16}(z)\subset B_{2R}(q),\qquad
 \Ric_g\ge-\varepsilon g,\qquad
 0\le\varepsilon\le1,\qquad \Vol(B_1(z))\ge v.
\]
Let \(u^1,\ldots,u^k,h\) be smooth real functions on \(B_8(z)\)
satisfying
\[
 \Delta u^a=0,\qquad \Delta h=2n,\qquad
 |u|+|du|+|h|+|dh|\le C_0.
\]
Set
\[
 \mathcal R=\Ric+\varepsilon g,\qquad
 b=h-|u|^2,\qquad \ell_a=\langle dh,du^a\rangle-2u^a,
\]
and define \(F\) to be the sum of the nonnegative integral
expressions on the left-hand sides of
\eqref{eq:common-energy} and \eqref{eq:homogeneity}, evaluated
with \(s=1\) and this \(\mathcal R\).
There are \(\eps_0=\eps_0(n,v,\eta,C_0)>0\) and
\(C=C(n,v,\eta,C_0)<\infty\) such that, if
\(\delta_0+F+\varepsilon\le\eps_0\), there is a constant
orthogonal matrix \(O\in O(k)\) such that,
after replacing \(u\) by \(Ou\), the functions
\[
 w'_j=u^{2j-1}+iu^{2j},\qquad 1\le j\le m-1,
\]
satisfy
\begin{equation}\label{eq:complex-pairing}
 \sum_{j=1}^{m-1}\int_{B_1(z)}|\bar\partial w'_j|^2\dd V_g
       \le C(F+\varepsilon).
\end{equation}
This change preserves \(b\) and \(\sum_a\ell_a^2\) exactly.
\end{lemma}

\begin{proof}
After choosing \(\eps_0\) sufficiently small, by the neck
symmetry gap and the almost splitting result, we have
\(\gamma=\gamma(n,v,\eta)>0\) such that \(B_{1/2}(z)\) admits
no harmonic \((k+1,\gamma)\)-splitting map. The conditions
\(r_z\le1/32\) and \(B_{16}(z)\subset B_{2R}(q)\) ensure
that the scales used here are above the stopping radius and
inside the neck's ambient ball.

Write \(X_a=\nabla u^a\) and
\(\widehat v_a=\frac12\langle JX_a,\nabla h\rangle\). Then
\[
 d\widehat v_a-(JX_a)^\flat
 =\tfrac12(\Hess h-2g)(JX_a,\cdot)
       +\tfrac12\langle J\Hess u^a,\nabla h\rangle,
 \qquad
 \int_{B_7(z)}|d\widehat v_a-(JX_a)^\flat|^2\dd V_g\le CF.
\]
Let \(v_a\) be its harmonic Dirichlet replacement on a smooth
domain between \(B_6(z)\) and \(B_7(z)\), and put
\(\xi_a=dv_a-(JX_a)^\flat\). Since \(\operatorname{div}(JX_a)=0\),
the Dirichlet energy estimate gives \(\int|\xi_a|^2\le CF\).
The action of \(J\) commutes with the Hodge Laplacian, so \(\xi_a\)
is Hodge harmonic. The Hodge Bochner and mean-value inequalities yield
\[
 \int_{B_3(z)}(|\xi_a|^2+|\nabla\xi_a|^2)\dd V_g
               +\sup_{B_2(z)}|\xi_a|^2\le CF.
\]
Define the constant skew matrix and residual forms by
\[
 \overline K_{ab}=\fintt_{B_2(z)}\langle JX_a,X_b\rangle\dd V_g,
 \qquad
 \beta_a=(JX_a)^\flat-\sum_b\overline K_{ab}\,du^b,\qquad
 A_a^2=\fintt_{B_2(z)}|\beta_a|^2\dd V_g.
\]
These forms are Hodge harmonic and \(A_a\le C\). Poincar\'e,
the squared Gram estimate, and the Hodge Bochner inequality give
\[
 \begin{split}
 &\int_{B_4(z)}|\nabla\beta_a|^2\dd V_g
 +\sum_b\int_{B_2(z)}|\langle\beta_a,du^b\rangle|^2\dd V_g
 +\fintt_{B_2(z)}(|\beta_a|-A_a)^2\dd V_g\le CF,\\
 &\sup_{B_1(z)}|\beta_a|\le CA_a.
 \end{split}
\]
For the cross terms use
\(\langle\beta_a,du^b\rangle
=\langle JX_a,X_b\rangle-\overline K_{ab}
-\sum_c\overline K_{ac}(\langle X_c,X_b\rangle-\delta_{cb})\).
For the length error, first use Poincar\'e to bound the variance.
The squared difference between the mean and the root mean square
is at most the variance.

If \(A_a>0\), set \(\tau=(F+\varepsilon)/A_a^2\).
Suppose \(\tau\) is sufficiently small. The harmonic function
\[
 w=\frac{v_a-\sum_b\overline K_{ab}u^b}{A_a},
 \qquad dw=\frac{\beta_a+\xi_a}{A_a}
\]
gives a harmonic map \(U=(u,w)\) satisfying
\[
 \int_{B_1(z)}
 \left[\sum_{i,j}|\langle dU^i,dU^j\rangle-\delta_{ij}|^2
                   +\sum_i|\Hess U^i|^2\right]\dd V_g\le C\tau,
 \qquad \sup_{B_1(z)}|dU^i|\le C.
\]
Here the normalized residual estimates follow by dividing by \(A_a\);
the old components contribute \(CF\le C\tau\).
The scalar Bochner inequality, \(\Delta|dU^i|^2\ge-C\varepsilon\),
and the mean-value inequality improve the gradient bound to
\(\sup_{B_{1/2}(z)}|dU^i|\le1+C\sqrt\tau\), since
\(\varepsilon\le A_a^2\tau\le C\tau\).
Thus \(U\) is a harmonic \((k+1,C\sqrt\tau)\)-splitting map,
contradicting the fixed gap. Consequently
\(\sum_a A_a^2\le C(F+\varepsilon)\).

Comparing the Gram matrices of \(JX_a\) and
\(\sum_b\overline K_{ab}X_b\) now gives
\[
 |\overline K\overline K^t-\mathrm{Id}_k|^2\le C(F+\varepsilon).
\]
Its polar part \(K_0=\overline K(-\overline K^2)^{-1/2}\)
is an orthogonal complex structure, with
\(|K_0-\overline K|^2\le C(F+\varepsilon)\).
A constant orthogonal change putting \(K_0\) into standard
complex blocks gives
\(\sum_j\int_{B_1(z)}|J\nabla u^{2j-1}-\nabla u^{2j}|^2
\dd V_g\le C(F+\varepsilon)\), which is
\eqref{eq:complex-pairing}. Orthogonality preserves \(|u|^2\)
and rotates the homogeneity error vector orthogonally.
\end{proof}

We return to the original metric and functions, keeping the
notation \(w'_j\) from Lemma~\ref{lem:complex-pairing}.
The following proposition combines the pairing estimate with
Corollary~\ref{cor:common-functions}; its energy bound controls
the errors and Ricci contractions in the proof of
Proposition~\ref{prop:regular}.

\begin{proposition}\label{prop:neck-energy}
Fix \(n=2m\ge4\) and \(v,\eta>0\), and put \(k=n-2\).
There is \(L_0=L_0(n,v,\eta)\ge32\) with the following
property. Fix \(L\ge L_0\), and let
\[
 \N=B_{2R}(q)\setminus\overline B_{r_x}(\C)
\]
be an \((n-2,\delta_0,\eta)\)-neck in a complete
K\"ahler manifold \((M^n,g,J)\) of real dimension \(n\).
Assume, for some \(\delta_1>0\),
\[
 \Ric_g\ge-(n-1)\delta_1^2g,\qquad
 \Vol(B_\rho(q))\ge v\rho^n\quad(0<\rho\le\delta_1^{-1}),
 \qquad \delta_1R\le4\delta_0\le1/100.
\]
Let \(\mu=\mu_{\N}\) be the packing measure of this neck,
and use \(\W^{\delta_1}\), defined by \eqref{eq:Wdefinition},
with one fixed cutoff based at \(q\). Take
\[
 y\in\N\cap B_R(q),\qquad z\in\C,\qquad
 r=d(y,\C)=d(y,z)>0,\qquad s=Lr\le10^{-6}R,
\]
and put \(\kappa=(n-1)\delta_1^2\). Let \(P\) be defined by
\eqref{eq:average-pinching}. There are
\(\eps_0>0\), \(C<\infty\), and
\(0<c\le C<\infty\), all depending only on
\(n,v,\eta,L\), such that, if
\(P+\kappa r^2+\delta_0\le\eps_0\), there exist smooth real
functions \(u^1,\ldots,u^k,h\) on \(B_{8s}(z)\) satisfying
\[
 \Delta u^a=0,\qquad \Delta h=2n,\qquad
 |u|\le Cs,\quad |du|\le C,\quad |h|\le Cs^2,\quad |dh|\le Cs,
\]
and
\begin{equation}\label{eq:lc-energy}
\begin{split}
 E:=r^{-n}\int_{B_r(y)}\bigg[
 &\sum_{a,b=1}^k|\langle du^a,du^b\rangle-\delta_{ab}|^2
       +\sum_{j=1}^{m-1}|\bar\partial w'_j|^2\\
 &+r^2\sum_{a=1}^k\Bigl(
     |\Hess u^a|^2+(\Ric+\kappa g)(\nabla u^a,\nabla u^a)\Bigr)\\
 &+r^{-2}\sum_{a=1}^k
       (\langle dh,du^a\rangle-2u^a)^2\\
 &+|\Hess h-2g|^2+(\Ric+\kappa g)(\nabla h,\nabla h)\\
 &+r^{-4}(|dh|^2-4h)^2\bigg]\dd V_g
       \le C(P+\kappa r^2).
\end{split}
\end{equation}

For \(b=h-|u|^2\), these same functions satisfy
\[
 cs^2\le b\le Cs^2\qquad\text{on }B_{r/8}(y).
\]
\end{proposition}

\begin{proof}
By Corollary~\ref{cor:common-functions}, with \(s=Lr\), we have
the stated functions, their pointwise bounds, and the bounds for
\(b\). It also bounds all terms in \eqref{eq:lc-energy} other
than the \(\bar\partial\)-term by \(CP\).
Divide lengths by \(s\), replace \(u,h\) by
\(s^{-1}u,s^{-2}h\), and apply
Lemma~\ref{lem:complex-pairing} in this normalization.
Indeed, the normalized stopping radius satisfies
\(r_z/s\le r/s=L^{-1}\le1/32\), and
\(B_{16s}(z)\subset B_{2R}(q)\).
Bishop--Gromov comparison and the volume assumption at \(q\)
give \(s^{-n}\Vol(B_s(z))\ge c(n)v\).
The normalized function bounds hold with
\(C_0=C_0(n,v,\eta,L)\), while
\[
 F\le CP,\qquad \varepsilon=\kappa s^2=L^2\kappa r^2.
\]
Thus the hypotheses of that lemma hold after decreasing
\(\eps_0(n,v,\eta,L)\). Restricting its conclusion to
\(B_r(y)\subset B_s(z)\) and returning to the original metric
gives
\[
 r^{-n}\int_{B_r(y)}\sum_{j=1}^{m-1}
    |\bar\partial w'_j|^2\dd V_g
 \le C\left(\frac{s}{r}\right)^n(F+\kappa s^2)
 \le C(P+\kappa r^2).
\]
The orthogonal change preserves \(|u|^2\) and all the other
sums, so it also preserves the stated bounds for \(b\).
\end{proof}

\section{A local Ricci estimate in terms of the energy}
\label{sec:regular}

This section proves a Ricci \(L^1\) estimate on \(B_r(y)\)
from functions and holomorphic coordinates on \(B_{8r}(y)\).
Proposition~\ref{prop:regular} gives a linear bound in their
energy \(E\). We fix all the local data below and use them
throughout this section.
Corollary~\ref{cor:neck-chart} and
Proposition~\ref{prop:neck-energy} provide these data on regular
balls in a neck; we make this application in Section~\ref{sec:proof}.

Let \(B_{8r}(y)\) be a relatively compact ball in a K\"ahler
manifold \((M^n,g,J)\) of real dimension \(n=2m\ge4\).
Fix \(C_0\ge1\) and \(\kappa\ge0\), and assume
\[
 \Ric_g\ge-\kappa g\quad\text{on }B_{8r}(y),\qquad
 \kappa r^2\le1.
\]
Suppose that there is a holomorphic chart
\(z:B_{8r}(y)\longrightarrow\mathbb C^m\), with \(z(y)=0\),
satisfying
\[
 \bigl||z(x)-z(x')|-d_g(x,x')\bigr|\le10^{-3}r
 \qquad(x,x'\in B_{8r}(y)),
\]
and
\[
 g\ge C_0^{-1}z^*g_{\mathbb C^m},\qquad
 \omega=i\partial\bar\partial\rho,\qquad |\rho|\le C_0r^2.
\]
We use the K\"ahler and coordinate conventions of
Section~\ref{sec:kahler-background}.
Let \(u^1,\ldots,u^{n-2},h\) be smooth real functions on
\(B_{8r}(y)\) such that
\[
 \Delta u^a=0,\qquad \Delta h=2n,\qquad b=h-|u|^2,
\]
and
\[
 C_0^{-1}r^2\le b\le C_0r^2,\qquad
 r^{-1}|u|+|du|+r^{-2}|h|+r^{-1}|dh|\le C_0.
\]
For the paired complex functions, we keep the notation \(w'_j\)
from Lemma~\ref{lem:complex-pairing}. Define
\begin{equation}\label{eq:local-energy}
\begin{split}
 E:=r^{-n}\int_{B_{8r}(y)}\bigg[
 &\sum_{a,b=1}^{n-2}|\langle du^a,du^b\rangle-\delta_{ab}|^2
       +\sum_{j=1}^{m-1}|\bar\partial w'_j|^2\\
 &+r^2\sum_{a=1}^{n-2}\Bigl(
     |\Hess u^a|^2+(\Ric+\kappa g)(\nabla u^a,\nabla u^a)\Bigr)\\
 &+r^{-2}\sum_{a=1}^{n-2}
       (\langle dh,du^a\rangle-2u^a)^2\\
 &+|\Hess h-2g|^2+(\Ric+\kappa g)(\nabla h,\nabla h)\\
 &+r^{-4}(|dh|^2-4h)^2\bigg]\dd V_g.
\end{split}
\end{equation}
In particular, \(E\) includes the nonnegative Ricci contractions
from the Bochner formulas, as well as all the squared errors
used below. Fix a smooth coordinate cutoff \(\chi\) with
\[
 \begin{gathered}
 0\le\chi\le1,\qquad \chi=1\text{ on }B_{2r}(y),\qquad
 \operatorname{supp}\chi\subset B_{3r}(y),\\
 |i\partial\bar\partial\chi|_0\le C(n)r^{-2}.
 \end{gathered}
\]
The coordinate distortion gives such a cutoff by rescaling a
fixed smooth radial profile.

The following proposition bounds the Ricci integral by this
local energy, with the error from the Ricci lower bound written
at scale \(r\).

\begin{proposition}\label{prop:regular}
In the local setting above, there is \(C=C(n,C_0)<\infty\)
such that
\begin{equation}\label{eq:local-bound}
 \int_{B_r(y)}|\Ric|\dd V_g
       \le Cr^{n-2}(E+\kappa r^2).
\end{equation}
\end{proposition}

\smallskip
\noindent\textit{Idea of the proof.}
We first explain the proof for \(E\le E_0(n,C_0)\).
The case of larger \(E\) follows from a coarse integral bound
using the same chart. Consider the section of the canonical bundle
\begin{equation}\label{eq:lc-sigma0}
 s_0=2^{-m/2}\partial w'_1\wedge\cdots\wedge
                 \partial w'_{m-1}\wedge\partial\log b.
\end{equation}
The product \(b|s_0|^2\) equals \(1\) on the model.
We estimate the two terms in
\[
 S+n\kappa
   =(S+n\kappa)b|s_0|^2
       +(S+n\kappa)(1-b|s_0|^2).
\]

\emph{The first term.}
The linear algebra in Lemma~\ref{lem:ricci-contractions} gives
\[
\begin{split}
 \int(S+n\kappa)b|s_0|^2\dd V_g
 &\le C\int\bigg[\sum_a
        (\Ric+\kappa g)(\nabla u^a,\nabla u^a)\\
 &\hspace{22mm}+r^{-2}
        (\Ric+\kappa g)(\nabla h,\nabla h)\bigg]\dd V_g
       \le Cr^{n-2}E.
\end{split}
\]
The radial term \(\Ric(\nabla h,\nabla h)\) appears in the
Bochner formula for \(h\). Its nonnegative shifted contraction,
and the corresponding contractions for \(u^a\), are included
in \(E\) with exactly the powers of \(r\) needed here.

\emph{The second term.}
This is the main difficulty, but its factor
\(1-b|s_0|^2\) has normalized \(L^2\) size \(\sqrt E\):
\[
 r^{-n}\int(1-b|s_0|^2)^2\dd V_g\le CE.
\]
Indeed, the squared norm is a Gram determinant, and \(E\) controls
the squared Gram and complex-pairing errors.
For this outline, we may assume \(s_0\) is holomorphic and
nonvanishing: Lemma~\ref{lem:holomorphic-correction} gives such
a replacement and bounds the curvature-weighted change by
\(Cr^{n-2}E\). We retain the name \(s_0\) for this replacement
only in the outline; the proof below denotes it by \(s\).

Positivity and the bounded K\"ahler potential first reduce the
scalar estimate to a coordinate cutoff integral, as in
\eqref{eq:lc-positive-comparison}.
Holomorphicity gives
\[
 2\pi\Ric(\omega)+\kappa\omega
   =i\partial\bar\partial\log(b|s_0|^2)
       -i\partial\bar\partial\log b+\kappa\omega.
\]
To estimate the remaining signed integral, we use the identity
\[
\begin{split}
 &(1-b|s_0|^2)\bigl(2\pi\Ric(\omega)+\kappa\omega\bigr)\\
 &\quad=(b|s_0|^2-1)i\partial\bar\partial\log b\\
 &\qquad+i b|s_0|^2\partial\log(b|s_0|^2)
                           \wedge\bar\partial\log(b|s_0|^2)\\
 &\qquad+\kappa(1-b|s_0|^2)\omega\\
 &\qquad-i\partial\bar\partial
                  \bigl[b|s_0|^2-1-\log(b|s_0|^2)\bigr].
\end{split}
\]
We integrate against \(\chi\omega_0^{m-1}\) in holomorphic
coordinates at scale \(r\), with a cutoff at the same scale.
The first term pairs two errors: the \((1,1)\)-Hessian
of \(\log b\) vanishes on the model, and its error is
linear in the Hessian and Gram errors. Thus
\[
 \|1-b|s_0|^2\|_{L^2}
       \|i\partial\bar\partial\log b\|_{L^2}
 \le C\bigl(r^{n/2}\sqrt E\bigr)
          \bigl(r^{n/2-2}\sqrt E\bigr)
 =Cr^{n-2}E.
\]
The second term is a squared covariant derivative. For the
original section, the product rule uses the cancellation
\[
 D'\partial\log b+\partial\log b\otimes\partial\log b
       =b^{-1}(\Hess b)^{2,0}.
\]
This vanishes on the model and is controlled by the Hessian
and complex-pairing errors. The holomorphic correction preserves
this estimate. The third term contains \(\kappa r^2\), and hence
\[
\begin{aligned}
 \int b|s_0|^2|\partial\log(b|s_0|^2)|^2\dd V_g
 &=\int b|D's_0+\partial\log b\otimes s_0|^2\dd V_g
       \le Cr^{n-2}E,\\
 \kappa\int|1-b|s_0|^2|\dd V_g
 &\le C\kappa r^n\sqrt E
       \le Cr^{n-2}(E+\kappa^2r^4).
\end{aligned}
\]
The metric lower bound gives the same bounds for the coordinate
cutoff integrals. For the last term, integrate by parts twice. The derivatives fall
on the cutoff, and the integral is with respect to \(dV_0\).
Writing \(s_0=2^{-m/2}a\,dz_1\wedge\cdots\wedge dz_m\),
\eqref{eq:lc-density} and Lemma~\ref{lem:holomorphic-correction} give
\[
 dV_0=\frac{b|s_0|^2}{b|a|^2}\,dV_g,
 \qquad b|a|^2\ge c.
\]
Thus the factor \(t=b|s_0|^2\) in
\(0\le t(t-1-\log t)\le(t-1)^2\) comes from the volume density.
The cutoff derivatives cost \(r^{-2}\), so
\[
\begin{split}
 Cr^{-2}\int\bigl[b|s_0|^2-1-\log(b|s_0|^2)\bigr]\dd V_0
 &\le Cr^{-2}\int(b|s_0|^2-1)^2\dd V_g\\
 &\le Cr^{n-2}E.
\end{split}
\]
Lemma~\ref{lem:local-calculation} estimates this remaining integral.
The final proof of Proposition~\ref{prop:regular} combines it
with Lemma~\ref{lem:ricci-contractions} and the positive-form
comparison to obtain, since \(\kappa r^2\le1\),
\[
 \int_{B_r(y)}|\Ric|\dd V_g
       \le Cr^{n-2}(E+\kappa r^2).
\]

\subsection{The weighted term}

We begin with the weighted scalar-curvature term, for which linear
algebra reduces the problem to Ricci contractions with
\(\nabla u^a\) and \(\nabla h\).
The right-hand side of the next lemma consists of the Ricci terms
included in the energy \(E\).

\begin{lemma}\label{lem:ricci-contractions}
In the local setting above, one has
\begin{equation}\label{eq:lc-sigma0-ricci}
\begin{split}
 &\int_{B_{8r}(y)}(S+n\kappa)b|s_0|^2\dd V_g\\
 &\quad\le C\int_{B_{8r}(y)}\bigg[
    \sum_{a=1}^{n-2}(\Ric+\kappa g)(\nabla u^a,\nabla u^a)
    +r^{-2}(\Ric+\kappa g)(\nabla h,\nabla h)
                                      \bigg]\dd V_g,
\end{split}
\end{equation}
where \(C=C(n,C_0)\).
\end{lemma}

\begin{proof}
Rescale \(g,u,h\) to \(r^{-2}g,r^{-1}u,r^{-2}h\), put
\(\varepsilon=\kappa r^2\), and retain the notation.
It suffices to prove the estimate for \(r=1\).

At a fixed point, choose an orthonormal eigenbasis
\(e_1,Je_1,\ldots,e_m,Je_m\) for \(\Ric+\varepsilon g\),
with eigenvalues \(r_1,r_1,\ldots,r_m,r_m\ge0\). Let \(\eta^1,\ldots,\eta^m\) be the unitary \((1,0)\)-coframe
dual to the vectors
\[
 \frac{e_\ell-iJe_\ell}{\sqrt2},\qquad 1\le\ell\le m.
\]
Write
\[
 \theta_j=2^{-1/2}\partial w'_j\quad(j<m),\qquad
 \theta_m=2^{-1/2}\partial\log b,\qquad
 \theta_j=\sum_\ell a_{j\ell}\eta^\ell.
\]
Then \(s_0=\theta_1\wedge\cdots\wedge\theta_m\).
For the Hermitian matrices
\[
 A_{\ell q}=\sum_j\overline{a_{j\ell}}a_{jq},\qquad
 R=\operatorname{diag}(r_1,\ldots,r_m),
\]
we have
\[
 \det A=|s_0|^2,\qquad
 \operatorname{tr}R=\tfrac12(S+n\varepsilon),\qquad
 \operatorname{tr}(AR)=\sum_{j,\ell}r_\ell|a_{j\ell}|^2.
\]
Since
\[
 d\log b=b^{-1}\left(dh-2\sum_a u^a\,du^a\right),
\]
the assumptions give \(0\le A\le C\mathrm{Id}_m\).
Diagonalizing \(A\), with eigenvalues \(\lambda_i\in[0,C]\),
and expressing \(R\) in the same basis gives
\begin{equation}\label{eq:lc-ricci-determinant}
\begin{split}
 \det A\operatorname{tr}R
 &=\sum_i\left(\prod_{j\ne i}\lambda_j\right)\lambda_iR_{i\bar i}
 \le C^{m-1}\sum_i\lambda_iR_{i\bar i}
 =C^{m-1}\operatorname{tr}(AR).
\end{split}
\end{equation}
This also holds when \(A\) is singular.

The unitary normalization gives the contraction identities
\[
\begin{split}
 \sum_\ell r_\ell|a_{j\ell}|^2
 &=\tfrac14(\Ric+\varepsilon g)
   (\nabla u^{2j-1}-J\nabla u^{2j},
                     \nabla u^{2j-1}-J\nabla u^{2j})\quad(j<m),\\
 \sum_\ell r_\ell|a_{m\ell}|^2
 &=\tfrac14(\Ric+\varepsilon g)(\nabla\log b,\nabla\log b).
\end{split}
\]
By positivity and \(J\)-invariance, we have
\[
 (S+n\varepsilon)|s_0|^2
 \le C\left[\sum_a(\Ric+\varepsilon g)(\nabla u^a,\nabla u^a)
       +(\Ric+\varepsilon g)(\nabla\log b,\nabla\log b)\right].
\]
Finally, by Cauchy--Schwarz for the nonnegative form
\(\Ric+\varepsilon g\), we have
\[
\begin{split}
 &(\Ric+\varepsilon g)(\nabla\log b,\nabla\log b)\\
 &\quad\le b^{-2}(1+4|u|^2)
   \left[(\Ric+\varepsilon g)(\nabla h,\nabla h)
                +\sum_a(\Ric+\varepsilon g)(\nabla u^a,\nabla u^a)\right].
\end{split}
\]
Multiplying the pointwise bound by \(b\le C\) and integrating
gives the estimate at unit scale. Restoring the original scale
gives \eqref{eq:lc-sigma0-ricci}.
\end{proof}

\begin{remark*}
The same estimate holds on any K\"ahler domain \(U\), with
\(U\) in place of \(B_{8r}(y)\), for arbitrary smooth
\(u:U\to\mathbb R^{n-2}\) and \(h:U\to\mathbb R\) satisfying
\[
 \begin{gathered}
 \Ric_g\ge-\kappa g,\qquad b=h-|u|^2,\qquad
 cr^2\le b\le C_0r^2,\\
 r^{-1}|u|+|du|+r^{-1}|dh|\le C_0,
 \end{gathered}
\]
where \(r>0\), \(\kappa\ge0\), and \(c>0\).
Here \(s_0\) is defined by \eqref{eq:lc-sigma0}, and the
constant depends only on \(n,c,C_0\). No harmonicity is required.
\end{remark*}

\subsection{The remaining term}

We now estimate the remaining term, whose small factor
\(1-b|s_0|^2\) must yield a linear bound in the energy \(E\).
Lemma~\ref{lem:local-calculation} gives the required signed integral
estimate, and Lemma~\ref{lem:holomorphic-correction} provides the
holomorphic replacement needed to prove it.

The next lemma controls the signed cutoff integral of the remaining
term at scale \(r\), with linear dependence on \(E\) and an explicit
error from the Ricci lower bound.

\begin{lemma}\label{lem:local-calculation}
Use the local setting and the cutoff fixed above.
There are \(E_0>0\) and \(C<\infty\), depending only on
\(n,C_0\), such that \(E\le E_0\) implies
\begin{equation}\label{eq:lc-conclusion}
\begin{split}
 \left|\int_{B_{3r}(y)}\chi(1-b|s_0|^2)
       \bigl(2\pi\Ric(\omega)+\kappa\omega\bigr)
                              \wedge\omega_0^{m-1}\right|
 \le Cr^{n-2}(E+\kappa^2r^4).
\end{split}
\end{equation}
\end{lemma}

To prove Lemma~\ref{lem:local-calculation}, we need a holomorphic
replacement for \(s_0\) with controlled energy and curvature-weighted
error.
The next lemma also gives nonvanishing on an interior ball, through
a lower bound for the coordinate coefficient that will be used to
compare \(dV_0\) and \(dV_g\).

\begin{lemma}\label{lem:holomorphic-correction}
Use the local setting fixed above.
There are \(E_0,c>0\) and \(C<\infty\), depending only on
\(n,C_0\), such that, if \(E\le E_0\), there is a
holomorphic section
\[
 s=2^{-m/2}a(z)\,dz_1\wedge\cdots\wedge dz_m
\]
of the canonical bundle on \(B_{6r}(y)\) such that
\begin{equation}\label{eq:lc-change-section}
 \int_{B_{4r}(y)}(S+n\kappa)b
                \bigl||s|^2-|s_0|^2\bigr|\dd V_g
       \le Cr^{n-2}E,
\end{equation}
\begin{equation}\label{eq:lc-j-energy}
 \int_{B_{4r}(y)}(b|s|^2-1)^2\dd V_g\le Cr^nE,
\end{equation}
and
\begin{equation}\label{eq:holomorphic-energy}
 \int_{B_{4r}(y)}\left(
       r^2|i\partial\bar\partial\log b|^2
       +b|D's+\partial\log b\otimes s|^2\right)\dd V_g
       \le Cr^{n-2}E.
\end{equation}
Moreover,
\begin{equation}\label{eq:lc-nonvanishing}
 r|a|\ge c\qquad\text{on }B_{4r}(y).
\end{equation}
\end{lemma}

\begin{proof}
Rescale \(g,u,h\) to \(r^{-2}g,r^{-1}u,r^{-2}h\), divide
the coordinates by \(r\), and retain the notation.
Put \(\varepsilon=\kappa r^2\le1\).
The section \(s_0\) rescales by \(r^{-(m-1)}\), so its norm
is multiplied by \(r\) and \(b|s_0|^2\) is unchanged.
The pointwise bounds in the local setting now hold with \(r=1\).
All constants in this proof depend only on \(n,C_0\).
The energy there is at most \(E\), volume comparison gives
\(\Vol(B_{8}(y))\le C(n)\), and the coordinate distortion gives
\[
 B_{6}(y)\subset B^0_{13/2}(0)
 \subset B^0_{7}(0)\subset B_{8}(y).
\]

Put \(\ell_a=\langle dh,du^a\rangle-2u^a\).
Expanding \(b=h-|u|^2\) gives
\begin{equation}\label{eq:lc-radial-identities}
\begin{split}
 |db|^2-4b
   ={}&(|dh|^2-4h)
     +4\sum_{a,b}u^au^b
          (\langle du^a,du^b\rangle-\delta_{ab})
                         -4\sum_a u^a\ell_a,\\
 \langle db,du^a\rangle
   ={}&\ell_a-2\sum_bu^b
                    (\langle du^a,du^b\rangle-\delta_{ab}).
\end{split}
\end{equation}
Together with the complex pairing term in \(E\),
these identities bound the
integral of the squared Gram error of
\(\nabla u^a,\nabla b/(2\sqrt b),J\nabla b/(2\sqrt b)\)
by \(CE\). Equality of the two Gram errors of a square matrix gives
\begin{equation}\label{eq:lc-frame}
 \int_{B_{8}(y)}\left|g-\sum_a du^a\otimes du^a
       -\frac{db\otimes db+(db\circ J)\otimes(db\circ J)}{4b}
       \right|^2\dd V_g\le CE.
\end{equation}
The Hessian identity
\begin{equation}\label{eq:lc-log-hessian}
\begin{split}
 &b\left(\Hess\log b+\tfrac12\left[
       d\log b\otimes d\log b
       -(d\log b\circ J)\otimes(d\log b\circ J)\right]\right)\\
 &\quad=(\Hess h-2g)-2\sum_a u^a\Hess u^a\\
 &\qquad\quad+2\left(g-\sum_a du^a\otimes du^a
       -\frac{db\otimes db+(db\circ J)\otimes(db\circ J)}{4b}\right)
\end{split}
\end{equation}
then gives, by taking its \((1,1)\) and \((2,0)\) parts,
\begin{equation}\label{eq:lc-log-energy}
 \int_{B_{8}(y)}\left(
    |i\partial\bar\partial\log b|^2
       +|D'\partial\log b
                    +\partial\log b\otimes\partial\log b|^2
                 \right)\dd V_g\le CE.
\end{equation}
Here the tensor in square brackets in
\eqref{eq:lc-log-hessian} is \(J\)-anti-invariant, and \(b\ge c\).

By the product rule in \eqref{eq:lc-sigma0}, the Hessian bounds
for \(u\), and \eqref{eq:lc-log-energy}, we have
\begin{equation}\label{eq:lc-sigma0-energy}
 \int_{B_{8}(y)}\left(|\bar\partial s_0|^2
       +|D's_0+\partial\log b\otimes s_0|^2\right)\dd V_g
       \le CE.
\end{equation}
The Gram determinant formula for \eqref{eq:lc-sigma0}, the
boundedness of its factors, and the preceding Gram estimates also give
\begin{equation}\label{eq:lc-j0}
 \int_{B_{8}(y)}(b|s_0|^2-1)^2\dd V_g\le CE.
\end{equation}

On the Stein coordinate ball \(B^0_{7}(0)\), apply the
\(L^2\) \(\bar\partial\) estimate in
\cite{Demailly2012} to the ordinary
\((m,1)\)-form \(\bar\partial s_0\), with trivial coefficient
bundle and weight \(e^{-\rho}\). Since
\(i\partial\bar\partial\rho=\omega\) and \(|\rho|\le C_0\),
we obtain
\begin{equation}\label{eq:lc-hormander}
\begin{split}
 &\bar\partial\xi=\bar\partial s_0,\\
 &\int_{B^0_{7}(0)}|\xi|^2\dd V_g
  \le e^{2C_0}
       \int_{B^0_{7}(0)}|\bar\partial s_0|^2\dd V_g
  \le CE.
\end{split}
\end{equation}
Thus \(s=s_0-\xi\) is holomorphic and \(\xi\) is smooth.

On the canonical bundle, the K\"ahler identities give
\begin{equation}\label{eq:lc-weitzenbock}
 2\bar\partial^*\bar\partial=D^*D+\frac S2;
\end{equation}
see \cite{Demailly2012}.
Test against \(\xi\) times the square of a coordinate cutoff
that equals \(1\) on \(B^0_{13/2}(0)\) and is supported in
\(B^0_{7}(0)\). The metric lower bound controls the cutoff
gradient; adding the scalar shift and absorbing cross terms gives
\begin{equation}\label{eq:lc-xi-energy}
\begin{split}
 &\int_{B_{6}(y)}\left(|D\xi|^2
                  +(S+n\varepsilon)|\xi|^2\right)\dd V_g\\
 &\quad\le C\int_{B^0_{7}(0)}
              \left(|\bar\partial s_0|^2+|\xi|^2\right)\dd V_g
       \le CE.
\end{split}
\end{equation}
The bound \(b\ge C_0^{-1}\), Lemma~\ref{lem:ricci-contractions},
\eqref{eq:lc-sigma0-energy}, and \(|d\log b|\le C\) imply
\begin{equation}\label{eq:lc-sigma-energy}
 \int_{B_{6}(y)}\left((S+n\varepsilon)|s|^2
        +|D's+\partial\log b\otimes s|^2\right)\dd V_g
       \le CE.
\end{equation}
Moreover, \(\bigl||s|^2-|s_0|^2\bigr|
\le2|s_0|\,|\xi|+|\xi|^2\), so weighted Cauchy--Schwarz gives
\[
\begin{split}
 &\int_{B_{6}(y)}(S+n\varepsilon)b
                  \bigl||s|^2-|s_0|^2\bigr|\dd V_g\\
 &\quad\le C
    \left(\int_{B_{6}(y)}
             (S+n\varepsilon)|s_0|^2\dd V_g\right)^{1/2}
    \left(\int_{B_{6}(y)}
             (S+n\varepsilon)|\xi|^2\dd V_g\right)^{1/2}\\
 &\qquad\quad+C\int_{B_{6}(y)}
                         (S+n\varepsilon)|\xi|^2\dd V_g
       \le CE.
\end{split}
\]

Writing \(s=2^{-m/2}a(z)\,dz_1\wedge\cdots\wedge dz_m\)
as in \cite{LiuSzekelyhidi2022}, we have
\begin{equation}\label{eq:lc-density}
 |s|^2=\frac{|a|^2}{\det(g_{a\bar b})},\qquad
                    |s|^2\dd V_g=|a|^2\dd V_0.
\end{equation}
The bounds for \(s_0\) and \(\xi\), holomorphic mean value
and Cauchy estimates, and the inverse metric bound give
\[
 \int_{B^0_{7}(0)}|a|^2\dd V_0\le C,\qquad
 \sup_{B_{6}(y)}
           (|a|+|da|_0+|s|+|\xi|)\le C.
\]
It follows from \eqref{eq:lc-j0} and \(s-s_0=-\xi\) that
\[
\begin{split}
 \int_{B_{6}(y)}(b|s|^2-1)^2\dd V_g
 &\le C\int_{B_{6}(y)}
             \left((b|s_0|^2-1)^2+|\xi|^2\right)\dd V_g
 \le CE.
\end{split}
\]

Choose \(\gamma>0\) with \(C_0^{m+1}\gamma^2\le1/4\).
If \(|a(z_*)|\le\gamma/2\) for some
\(z_*\in B_{4}(y)\), the Cauchy bound gives a coordinate
ball \(B^0_{\rho_*}(z_*)\subset B_{6}(y)\), with
\(\rho_*=\rho_*(n,C_0)>0\), on which \(|a|\le\gamma\).
There \(b|s|^2\le1/4\) and \(dV_g\ge C_0^{-m}dV_0\), so
\[
 \int_{B_{6}(y)}(b|s|^2-1)^2\dd V_g
 \ge\frac9{16}C_0^{-m}
                         \Vol_0(B^0_{\rho_*})>0.
\]
This contradicts the preceding bound for
\(E\le E_0(n,C_0)\) small enough.
Thus \(|a|\ge\gamma/2\) on \(B_{4}(y)\).

Restoring the original metric multiplies \(|s_0|,|s|,|a|\)
by \(r^{-1}\) and \(b\) by \(r^2\). This gives the powers
of \(r\) in \eqref{eq:lc-change-section}--\eqref{eq:holomorphic-energy}
and the coefficient bound \eqref{eq:lc-nonvanishing}.
\end{proof}

\begin{proof}[Proof of Lemma~\ref{lem:local-calculation}]
Rescale \(g,u,h\) to \(r^{-2}g,r^{-1}u,r^{-2}h\), divide
the holomorphic coordinates by \(r\), and retain the same
notation. Put \(\varepsilon=\kappa r^2\le1\). The quantity
\(E\) is unchanged, and it suffices to prove a bound
\(C(E+\varepsilon^2)\) at this scale.

Apply Lemma~\ref{lem:holomorphic-correction} to obtain the
nonvanishing holomorphic section \(s\) on
\(B_{4}(y)\). Since
\(\omega_0\le C\omega\) and
\(2\pi\Ric(\omega)+\varepsilon\omega\ge0\),
\eqref{eq:lc-change-section} gives
\[
\begin{split}
 &\left|\int\chi b\bigl(|s|^2-|s_0|^2\bigr)
       \bigl(2\pi\Ric(\omega)+\varepsilon\omega\bigr)
                                      \wedge\omega_0^{m-1}\right|\\
 &\quad\le C\int_{B_{4}(y)}
       (S+n\varepsilon)b\bigl||s|^2-|s_0|^2\bigr|\dd V_g
       \le CE.
\end{split}
\]
Thus it remains to estimate the cutoff integral with \(s\)
in place of \(s_0\).

Write \(s=2^{-m/2}a\,dz_1\wedge\cdots\wedge dz_m\).
Since \(a\) is holomorphic and nonvanishing,
\eqref{eq:lc-density} gives
\[
 i\partial\bar\partial\log|s|^2=2\pi\Ric(\omega).
\]
Consequently,
\[
\begin{split}
 i\partial\bar\partial\log(b|s|^2)
   &=2\pi\Ric(\omega)+i\partial\bar\partial\log b,\\
 i\partial\bar\partial(b|s|^2)
   &=b|s|^2 i\partial\bar\partial\log(b|s|^2)\\
   &\quad+i b|s|^2\partial\log(b|s|^2)
                         \wedge\bar\partial\log(b|s|^2).
\end{split}
\]
Subtracting the second identity from the first and adding
\(\varepsilon(1-b|s|^2)\omega\) gives
\begin{equation}\label{eq:lc-log-identity}
\begin{split}
 &(1-b|s|^2)\bigl(2\pi\Ric(\omega)+\varepsilon\omega\bigr)\\
 &\quad=(b|s|^2-1)i\partial\bar\partial\log b\\
 &\qquad+i b|s|^2\partial\log(b|s|^2)
                              \wedge\bar\partial\log(b|s|^2)\\
 &\qquad+\varepsilon(1-b|s|^2)\omega\\
 &\qquad-i\partial\bar\partial
                        \bigl[b|s|^2-1-\log(b|s|^2)\bigr].
\end{split}
\end{equation}
We integrate this identity against \(\chi\omega_0^{m-1}\).

For the first term, by \eqref{eq:lc-j-energy} and
\eqref{eq:holomorphic-energy}, we have
\[
\begin{split}
 &\int_{B_{3}(y)}\bigl|b|s|^2-1\bigr|
                           |i\partial\bar\partial\log b|\dd V_g\\
 &\quad\le
  \left(\int_{B_{4}(y)}(b|s|^2-1)^2\dd V_g\right)^{1/2}
  \left(\int_{B_{4}(y)}
                 |i\partial\bar\partial\log b|^2\dd V_g\right)^{1/2}
  \le CE.
\end{split}
\]
The bound \(\omega_0\le C\omega\) controls the absolute value
of its cutoff integral by this expression.

For the second term, the Chern connection formula
\cite{Demailly2012} gives
\begin{equation}\label{eq:lc-connection}
\begin{split}
 D's+\partial\log b\otimes s
       &=\partial\log(b|s|^2)\otimes s,\\
 b|s|^2|\partial\log(b|s|^2)|^2
       &=b|D's+\partial\log b\otimes s|^2.
\end{split}
\end{equation}
Its nonnegative cutoff integral is therefore bounded by
\[
\begin{split}
 C\int_{B_{3}(y)}
             b|s|^2|\partial\log(b|s|^2)|^2\dd V_g
 &=C\int_{B_{3}(y)}
             b|D's+\partial\log b\otimes s|^2\dd V_g\\
 &\le CE
\end{split}
\]
by \eqref{eq:holomorphic-energy}. For the third term, the
volume bound and \eqref{eq:lc-j-energy} imply
\[
 \varepsilon\int_{B_{3}(y)}
              \bigl|1-b|s|^2\bigr|\dd V_g
 \le C\varepsilon\sqrt E
 \le C(E+\varepsilon^2).
\]

For the last term, use
\[
 0\le t(t-1-\log t)\le(t-1)^2\qquad(t>0),
\]
which follows from
\((t-1)^2-t(t-1-\log t)=t\log t-t+1\ge0\).
The density identity \eqref{eq:lc-density},
\eqref{eq:lc-nonvanishing}, and \(b\ge c\) give
\begin{equation}\label{eq:lc-remainder}
\begin{split}
 &\bigl[b|s|^2-1-\log(b|s|^2)\bigr]\dd V_0\\
 &\quad=\frac{b|s|^2\bigl[b|s|^2-1-\log(b|s|^2)\bigr]}
                    {b|a|^2}\dd V_g
 \le C(b|s|^2-1)^2\dd V_g.
\end{split}
\end{equation}
In particular, the first factor \(t=b|s|^2\) in the
elementary inequality comes from the volume density.
Integration by parts in the coordinates gives
\[
\begin{split}
 &\left|\int\chi\,i\partial\bar\partial
       \bigl[b|s|^2-1-\log(b|s|^2)\bigr]\wedge\omega_0^{m-1}\right|\\
 &\quad=\left|\int\bigl[b|s|^2-1-\log(b|s|^2)\bigr]
                    i\partial\bar\partial\chi\wedge\omega_0^{m-1}\right|\\
 &\quad\le C\int_{B_{3}(y)}
             \bigl[b|s|^2-1-\log(b|s|^2)\bigr]\dd V_0\\
 &\quad\le C\int_{B_{3}(y)}(b|s|^2-1)^2\dd V_g\le CE.
\end{split}
\]
Here the coordinate Hessian bound for \(\chi\) controls
the measure in the second line by \(C\dd V_0\).
No pointwise lower bound for \(b|s|^2\) is needed.

Combining these four estimates in
\eqref{eq:lc-log-identity}, and then using the replacement
estimate at the beginning of the proof, gives
\[
 \left|\int\chi(1-b|s_0|^2)
       \bigl(2\pi\Ric(\omega)+\varepsilon\omega\bigr)
                                      \wedge\omega_0^{m-1}\right|
       \le C(E+\varepsilon^2).
\]
Restoring the original scale multiplies this integral by
\(r^{n-2}\). Since \(\varepsilon=\kappa r^2\), this proves
\eqref{eq:lc-conclusion}.
\end{proof}

\subsection{Completion of the local estimate}

We combine the two local estimates to prove
Proposition~\ref{prop:regular}.
Positivity and the bounded K\"ahler potential compare the
scalar integral with the coordinate cutoff integrals.
The same comparison also gives the coarse bound needed when
\(E\) is not small.

\begin{proof}[Proof of Proposition~\ref{prop:regular}]
The coordinate distortion and a first-exit argument give
\[
 B_r(y)\subset B^0_{5r/4}(0)
       \subset B^0_{3r/2}(0)\subset B_{2r}(y).
\]
The Ricci form is nonnegative after the shift, and
\[
 \bigl(2\pi\Ric(\omega)+\kappa\omega\bigr)\wedge\omega^{m-1}
       =\frac{S+n\kappa}{2m}\omega^m.
\]
Apply the Chern--Levine--Nirenberg inequality
\cite{Demailly2012} on the displayed
coordinate balls. The bound \(|\rho|\le C_0r^2\) cancels
the factors \(r^{-2}\) from the cutoff Hessians, so
\begin{equation}\label{eq:lc-positive-comparison}
\begin{split}
 \int_{B_r(y)}(S+n\kappa)\dd V_g
 &\le C\int_{B^0_{3r/2}(0)}
       \bigl(2\pi\Ric(\omega)+\kappa\omega\bigr)
                                      \wedge\omega_0^{m-1}\\
 &\le C\int\chi\bigl(2\pi\Ric(\omega)+\kappa\omega\bigr)
                                      \wedge\omega_0^{m-1}.
\end{split}
\end{equation}
All constants below depend only on \(n,C_0\).

Suppose first that \(E\le E_0(n,C_0)\), so that
Lemmas~\ref{lem:local-calculation} and~\ref{lem:holomorphic-correction}
apply. For the weighted term, \(\omega_0\le C\omega\) and
Lemma~\ref{lem:ricci-contractions} give
\[
\begin{split}
 &\int\chi b|s_0|^2
       \bigl(2\pi\Ric(\omega)+\kappa\omega\bigr)
                                      \wedge\omega_0^{m-1}\\
 &\quad\le C\int_{B_{8r}(y)}\bigg[
       \sum_a(\Ric+\kappa g)(\nabla u^a,\nabla u^a)
       +r^{-2}(\Ric+\kappa g)(\nabla h,\nabla h)
                             \bigg]\dd V_g
       \le Cr^{n-2}E.
\end{split}
\]
Use Lemma~\ref{lem:local-calculation} for the remaining term.
Together with \(|\Ric|\le S+2n\kappa\),
\eqref{eq:lc-positive-comparison}, and volume comparison,
this yields
\begin{equation}\label{eq:local-rescaled-final}
\begin{split}
 \int_{B_r(y)}|\Ric|\dd V_g
 &\le C\int\chi b|s_0|^2
       \bigl(2\pi\Ric(\omega)+\kappa\omega\bigr)
                                      \wedge\omega_0^{m-1}\\
 &\quad+C\left|\int\chi(1-b|s_0|^2)
       \bigl(2\pi\Ric(\omega)+\kappa\omega\bigr)
                            \wedge\omega_0^{m-1}\right|+C\kappa r^n\\
 &\le Cr^{n-2}(E+\kappa^2r^4)+C\kappa r^n\\
 &\le Cr^{n-2}(E+\kappa r^2),
\end{split}
\end{equation}
since \(\kappa r^2\le1\).

For completeness, the chart gives a coarse bound without any
smallness of \(E\). Put \(G=(g_{a\bar b})\) and
\(f=\log\det G\) in the holomorphic coordinates. By the metric
lower bound, we have
\[
 |f|\le\det G+m\log C_0,\qquad
 \operatorname{tr}G\le mC_0^{m-1}\det G.
\]
Since \(dV_g=\det G\,dV_0\), the coordinate distortion and
volume comparison imply
\[
 \int_{B_{3r}(y)}\bigl(|f|+\operatorname{tr}G\bigr)\dd V_0
       \le Cr^n.
\]
Using \(2\pi\Ric(\omega)=-i\partial\bar\partial f\),
integration by parts gives
\[
\begin{split}
 \int\chi\bigl(2\pi\Ric(\omega)+\kappa\omega\bigr)
                  \wedge\omega_0^{m-1}
 &=-\int f\,i\partial\bar\partial\chi\wedge\omega_0^{m-1}
       +\kappa\int\chi\omega\wedge\omega_0^{m-1}\\
 &\le Cr^{-2}\int_{B_{3r}(y)}|f|\dd V_0
       +C\kappa\int_{B_{3r}(y)}\operatorname{tr}G\dd V_0\\
 &\le Cr^{n-2}.
\end{split}
\]
Thus \eqref{eq:lc-positive-comparison} gives
\(\int_{B_r(y)}|\Ric|\dd V_g\le Cr^{n-2}\).
When \(E>E_0\), this is at most
\(CE_0^{-1}r^{n-2}E\), proving \eqref{eq:local-bound}
in this case as well.
\end{proof}

\section{Proof of the main theorem}
\label{sec:proof}

We now complete the proof of Theorem~\ref{thm:main}, following
the neck summation argument of Jiang--Naber~\cite{JiangNaber2021}.
The new input is the local Ricci $L^1$ estimate from
Section~\ref{sec:regular}. Combined with
Proposition~\ref{prop:neck-energy}, it gives the required entropy
estimate on regular balls, which can then be summed by the neck
structure and entropy monotonicity.

Fix \(\eta=\eta(n,v)>0\), then choose
\(a_0,L,\delta_0\) as in the preceding results.
In this section the manifolds are K\"ahler of real dimension \(n=2m\ge4\).

The next proposition gives the \(L^1\) bound for a whole neck.
Its proof uses the Ahlfors estimate in
Theorem~\ref{thm:neck-structure} and telescoping of the entropy
drops, with Theorem~\ref{thm:scalar-almost-euclidean} handling the
remaining almost Euclidean balls.

\begin{proposition}\label{prop:neck}
For \(n=2m\ge4\) and \(v>0\), there are \(\eta=\eta(n,v)>0\) and \(\delta=\delta(n,v)>0\)
with the following property. Let
\(\N=B_{2R}(q)\setminus\overline B_{r_x}(\C)\) be an \((n-2,\delta_0,\eta)\)-neck in a complete
K\"ahler manifold \((M^n,g,J)\), with
\(0<\delta_0\le\delta\). Suppose, for \(\delta_1>0\), that
\[
 \Ric_g\ge-(n-1)\delta_1^2g,\qquad
 \Vol(B_\rho(q))\ge v\rho^n\quad(0<\rho\le\delta_1^{-1}),
 \qquad \delta_1R\le4\delta_0\le1/100.
\]
Then
\begin{equation}\label{eq:neck-bound}
 \int_{\N\cap B_R(q)}|\Ric|\dd V_g\le C(n,v)R^{n-2}.
\end{equation}
\end{proposition}

\begin{proof}
Fix \(L=L(n,v,\eta)\ge L_0(n,v,\eta)\) as in
Proposition~\ref{prop:neck-energy}, and choose
\(a_0=a_0(n,v,\eta)\) and the chart accuracy as after
Corollary~\ref{cor:neck-chart}. Put
\(\kappa=(n-1)\delta_1^2\) and \(\mu=\mu_{\N}\).
All entropy values below use one fixed cutoff based at \(q\).

We first obtain the estimate on each regular ball. For
\(y\in\N\cap B_R(q)\), choose \(z\in\C\) such that
\begin{equation}\label{eq:regular-radius}
 r=d(y,\C)=d(y,z)>0,\qquad
 r\le\frac{R}{10^6L},\qquad s=Lr,
\end{equation}
and suppose that \(B_{a_0r}(y)\subset\N\).
Take \(P\) from \eqref{eq:average-pinching}.
Choose a sufficiently small \(c=c(n,v,\eta,L)>0\), and decrease
\(\delta(n,v)\) so that Proposition~\ref{prop:neck-energy}
applies whenever \(P\le c\). This is possible since
\(\kappa r^2\le16(n-1)\delta_0^2\).

For \(P\le c\), restrict the functions of
Proposition~\ref{prop:neck-energy} to \(B_{8a_0r}(y)\).
The chart in Corollary~\ref{cor:neck-chart}, written in the
original metric, and the bounds for these functions give the
local setting of Section~\ref{sec:regular} with radius \(a_0r\)
and \(C_0=C_0(n,v,\eta,L)\). Here
\(8a_0<1/8\), so the positive lower bound for \(b\) holds
on the whole chart. Each term of the local energy
\eqref{eq:local-energy}, evaluated at radius \(a_0r\),
is at most \(a_0^{-n-4}\) times its counterpart in
\eqref{eq:lc-energy}. By Proposition~\ref{prop:regular}, we have
\[
 \int_{B_{a_0r}(y)}|\Ric|\dd V_g
       \le Cr^{n-2}(P+\kappa r^2).
\]
If \(P>c\), Theorem~\ref{thm:scalar-almost-euclidean}
on the same regular ball gives
\[
 \int_{B_{a_0r}(y)}|\Ric|\dd V_g
       \le Cr^{n-2}\le Cr^{n-2}P.
\]
The chart accuracy gives the required Euclidean approximation;
volume comparison and a further decrease of \(\delta\)
give the volume and Ricci hypotheses.
Since \(r^{n-2}P=L^{2-n}s^{n-2}P\), both cases yield
\begin{equation}\label{eq:H}
\begin{split}
 \int_{B_{a_0r}(y)}|\Ric|\dd V_g
 \le{}&A\int_{\C\cap B_s(z)}
   \left(\W_{s^2/2}^{\delta_1}(x)
                    -\W_{2s^2}^{\delta_1}(x)\right)\dd\mu(x)
       +A\delta_1^2r^n,
\end{split}
\end{equation}
where \(A=A(n,v,\eta,L)\). Since \(\eta=\eta(n,v)\)
and \(L=L(n,v,\eta)\) are fixed, all constants depend only
on \(n,v\).

Put \(\theta=(10^6L)^{-1}\). Choose balls
\(B_{a_0r_i}(y_i)\), with disjoint
\(B_{a_0r_i/5}(y_i)\), to cover the points with
\(d(y,\C)\le\theta R\) and \(d(y,x)>100r_x\) for every
\(x\in\C\). Here
\(r_i=d(y_i,\C)=d(y_i,z_i)\) and \(s_i=Lr_i\).
The factor \(100\) ensures that these balls lie in \(\N\).
Volume comparison gives \(\sum_i r_i^n\le CR^n\) and bounded
overlap, at each dyadic scale, of the center sets
\(\C\cap B_{s_i}(z_i)\). Their union has packing measure at
most \(CR^{n-2}\) by Theorem~\ref{thm:neck-structure}; all its
centers lie in \(B_{3R/2}(q)\) and have \(r_x\le R/16\).

For \(\rho_j=2^{-j}\theta R\), \(T_j=(L\rho_j)^2\), and
\(\rho_j/2<r_i\le\rho_j\), monotonicity gives
\[
 0\le\W_{s_i^2/2}^{\delta_1}(x)-\W_{2s_i^2}^{\delta_1}(x)
 \le\W_{T_j/8}^{\delta_1}(x)-\W_{2T_j}^{\delta_1}(x).
\]
Since \(2T_{j+3}<T_j/8\), we sum separately over the three
residue classes of \(j\) modulo three. Each sum telescopes
at a fixed center. By Remark~\ref{rem:entropy-bound}
and \eqref{eq:H}, we have
\begin{equation}\label{eq:main-sum}
\begin{split}
 \sum_i\int_{B_{a_0r_i}(y_i)}|\Ric|\dd V_g
 &\le A\sum_i\int_{\C\cap B_{s_i}(z_i)}
       \left(\W_{s_i^2/2}^{\delta_1}(x)
                   -\W_{2s_i^2}^{\delta_1}(x)\right)\dd\mu(x)
          +CA\delta_1^2R^n\\
 &\le CR^{n-2}.
\end{split}
\end{equation}
The remaining points have \(d(y,\C)>\theta R\), or lie in
stopping collars where \(r_x/2\le d(y,\C)\le100r_x\).
Corollary~\ref{cor:neck-chart} covers them by almost Euclidean
balls at fixed multiples of \(\theta R\) and \(r_x\), with
total \((n-2)\)-content at most \(CR^{n-2}\) by the same
Ahlfors estimate. Theorem~\ref{thm:scalar-almost-euclidean} completes
\eqref{eq:neck-bound}.
\end{proof}

The bounds for necks and almost Euclidean balls now have the same
\(r^{n-2}\) scaling, so the packing estimate in
Theorem~\ref{thm:neck-decomposition} completes the proof.

\begin{proof}[Proof of Theorem~\ref{thm:main}]
By choosing
\(\eta=\eta(n,v)>0\) sufficiently small, the codimensional two singularity implies that every
\((n-1,2\eta)\)-symmetric ball in the decomposition is indeed almost
regular. Hence, Theorem~\ref{thm:scalar-almost-euclidean} applies to these
balls. Proposition~\ref{prop:neck} bounds each neck. By Theorem~\ref{thm:neck-decomposition}, the estimate 
\eqref{eq:packing-content} gives
\begin{equation}\label{eq:global-sum}
\begin{split}
 \int_{B_1(p)}|\Ric|\dd V_g
 &\le\sum_a\int_{\N_a\cap B_{r_a}(x_a)}|\Ric|\dd V_g
       +\sum_b\int_{B_{r_b}(x_b)}|\Ric|\dd V_g\\
 &\le C\left(\sum_a r_a^{n-2}+\sum_b r_b^{n-2}\right)
 \le C(n,v).
\end{split}
\end{equation}
\end{proof}

\section*{Acknowledgments}
\enlargethispage{3\baselineskip}

The authors are deeply grateful to Aaron Naber for his valuable
suggestions and many discussions on this problem over the years.
His ideas and insights have been an important source of inspiration
for this work and have helped shape our approach to the problem.
We also thank him for his encouragement and generosity in sharing
his mathematical ideas.

The authors acknowledge the substantial assistance of OpenAI's
ChatGPT in developing the local linear estimate in
Proposition~\ref{prop:regular}. In earlier attempts, the authors
had obtained only a square-root estimate in the entropy error.
AI-assisted exploration led to the decomposition of the
scalar-curvature integral described in the ``Idea of the proof''
of Proposition~\ref{prop:regular}. This decomposition was a key
step in obtaining the linear bound required for the neck summation.
The authors subsequently checked and completed the argument and
reorganized its presentation. They take full responsibility for
the mathematical content of the paper.

This work was supported by the National Key Research and Development Program
of China (No.~2022YFA1005500), the National Natural Science Foundation of
China (Grant Nos.~12125105 and 12071140), the New Cornerstone Science
Foundation through the New Cornerstone Investigator Program, the Xplorer
Prize, and the Science and Technology Commission of Shanghai Municipality
(No.~22DZ2229014).

\begingroup
\small\linespread{0.96}\selectfont

\endgroup
\end{document}